\documentclass[a4paper,11pt]{article}
\usepackage{amsfonts,amssymb,latexsym,amsmath}
\usepackage{multicol}
\newcommand{\GL}{{\mathrm{GL}(n,K)}}
\newcommand{\Ad}{{\mathrm{Ad}}}
\newcommand{\Lc}{{\cal L}}
\newcommand{\Xb}{{\Bbb X}}
\newcommand{\Yb}{{\Bbb Y}}
\newcommand{\Sb}{{\Bbb S}}
\renewcommand{\leq}{\leqslant}

\begin{document}
\Large

\title{Field of U-invariants of adjoint representation of the group
$\GL$}
\author{K.A.Vyatkina \and A.N.Panov \thanks{The work is  supported by the
RFBR-grants  12-01-00070-a, 12-01-00137-a, 13-01-97000-Volga
region-a}}

\date{}

\maketitle

It is known that the field of invariants  of an arbitrary
unitriangular group is rational  (see \cite{m}). In this paper for
adjoint representation of the group $\GL$ we present the
algebraically system of generators of the field of $U$-invariants.
Note that the algorithm  for calculating  generators of the field of
invariants of an arbitrary rational representation of an unipotent
group was presented in  \cite[Chapter 1]{p}. The invariants was
constructed using induction on the length of Jordan-H\"older series.
In our case the length is equal to  $n^2$; that is why it is
difficult to apply the  algorithm.

 Let us consider the adjoint representation  $\mathrm{Ad}_gA =gAg^{-1}$ of the group
 $\mathrm{GL}(n,K)$, where $K$ is a field of zero characteristic, on the algebra of matrices $M=\mathrm{Mat}(n,K)$.
 The adjoint representation determines the  representation  $\rho_g$  on   $K[M]$ (resp. $K(M)$)
 by formula
$\rho_gf(A) = f(g^{-1} A g).$

Let $U$ be the subgroup of upper triangular matrices in
$\mathrm{GL}(n,K)$ with units on the diagonal. A polynomial
(rational function) $f$ on $M$ is called an $U$-invariant, if
$\rho_u f = f$ for every $u\in U$. The  set of $U$-invariant
rational functions $K(M)^U$ is a subfield of  $K(M)$.

Let  $\{x_{i,j}\}$  be a system of standard coordinate functions on
 $M$. Construct a matrix  $\Xb = (x_{ij})$. Let $\Xb^* = (x^*_{ij})$ be its adjugate
matrix, $ \Xb\cdot \Xb^*= \Xb^*\cdot \Xb =\det \Xb\cdot E$. Denote
by  $J_k$ the left lower corner minor of order $k$ of the matrix
$\Xb$.

 We associate with each  $J_k$ the system of $k$ determinants  $J_{k,i}$, where
$0\leq i\leq k-1$. The determinant $J_{k,0}$ coincides with the
minor $J_k$. The first $k-i$ rows in the determinant  $J_{k,i}$,
coincide with the last  $k-i$ rows in the minor $J_k$; the last $i$
rows in $J_{k,i}$ coincide with the similar rows in the minor
$J_k(\Xb^*)$ of the adjugate matrix  $\Xb^* = (x_{ij}^*)$. \\
\\
EXAMPLE  1. Case $ n=2$,\quad  $\Xb=\left(\begin{array}{cc}
x_{11}&x_{12}\\
x_{21}&x_{22}\end{array}\right)$,\quad
$\Xb^*=\left(\begin{array}{cc}
x^*_{11}&x^*_{12}\\
x^*_{21}&x^*_{22}\end{array}\right)$\\
\\
$ J_{1,0} = x_{2,1}$,\,\, $J_{2,0} = \left|\begin{array}{cc}
x_{11}&x_{12}\\
x_{21}&x_{22}\end{array}\right|$,\quad $J_{2,1} =
\left|\begin{array}{cc}
x_{21}&x_{22}\\
x^*_{21}&x^*_{22}\end{array}\right|$.\\
\\
\\
EXAMPLE  2. Case $ n=3$, \\
\\
$\Xb=\left(\begin{array}{ccc}
x_{11}&x_{12}&x_{13}\\
x_{21}&x_{22}& x_{23}\\
x_{31}&x_{32}& x_{33}\end{array}\right)$, \quad
$\Xb^*=\left(\begin{array}{ccc}
x^*_{11}&x^*_{12}&x^*_{13}\\
x^*_{21}&x^*_{22}& x^*_{23}\\
x^*_{31}&x^*_{32}& x^*_{33}\end{array}\right)$,
\\
\\
\\
$ J_{1,0} = x_{31}$,\quad $J_{2,0} = \left|\begin{array}{cc}
x_{21}&x_{22}\\
x_{31}&x_{32}\end{array}\right|$,\quad $J_{3,0} =
\left|\begin{array}{ccc}
x_{11}&x_{12}&x_{13}\\
x_{21}&x_{22}& x_{23}\\
x_{31}&x_{32}& x_{33}\end{array}\right|$,\\
\\
\\
$J_{2,1} = \left|\begin{array}{cc}
x_{31}&x_{32}\\
x^*_{31}&x^*_{32}\end{array}\right|$,\quad $J_{3,1} =
\left|\begin{array}{ccc}
x_{21}&x_{22}&x_{23}\\
x_{31}&x_{32}& x_{33}\\
x^*_{31}&x^*_{32}& x^*_{33}\end{array}\right|$,\quad $J_{3,2} =
\left|\begin{array}{ccc}
x_{31}&x_{32}&x_{33}\\
x^*_{21}&x^*_{22}& x^*_{23}\\
x^*_{31}&x^*_{32}& x^*_{33}\end{array}\right|$.\\
\\
\\
{\bf Proposition 1}. Each polynomial  $J_{k,i}$, where  $ 1\leq
k\leq n$,~ $0\leq i\leq k-1$ is an   $U$-invariant.
\\
{\bf Proof}. First, the formula $\rho_u\Xb = u^{-1}\Xb u$ implies
that each left lower corner minor  $J_k$, where $1\leq k\leq n$, is
an  $U$-invariant. Indeed, each $J_k$ is invariant with respect to
right and left  multiplication by elements  $u\in U$. For each
$1\leq i\leq n$ we construct the  $n\times n $ matrix
$$\Yb_i =
\left(\begin{array}{c} \Xb_{n-i}\\ ----- \\
\Xb^*_{i}\end{array}\right).$$ Here  the block $\Xb_{n-i}$ has size
$(n-i)\times n$; it is obtained from the matrix $\Xb$ deleting first
$i$ rows. The block $\Xb^*_{i}$ has size $i\times n$; it is obtained
from the matrix $\Xb^*$ by deleting first $n-i$ rows.

Since  $\rho_u\Xb = u^{-1}\Xb u$ and $\rho_u\Xb^* = u^{-1}\Xb^* u$,
we have
$$\rho_u \Yb_i =
\left(\begin{array}{c} u_{n-i}^{-1}\Xb_{n-i}u\\ ----- \\
u_{i}^{-1}\Xb^*_{i}u\end{array}\right) = u_*\Yb_iu,\eqno{(1)}$$
 where $ u_* = \left(\begin{array}{cc}
u_{n-i}^{-1}& 0\\
0& u_i^{-1}\end{array}\right)$, ~~  $u_{n-i}$ (resp. $u_{i}$) is a
right lower block of order  $n-i$ (resp. $i$) for $u\in U$.

Note that for  every $1\leq i\leq k-1$ the determinant  $J_{k,i}$ is
the left lower corner minor of order $k$ of the matrix  $\Yb_i$. The
formula (1) implies that $J_{k,i}$ is an $U$-invariant. $\Box$

Consider a Zariski open subset  $\Omega\subset M$, that consists of
all matrices obeying  $J_k\ne 0$ for any $1\leq k\leq n$. Denote
 by  $\Lc$ the subset  of matrices of form  $$ B =
\left(\begin{array}{ccccc}
0&0&0&\ldots&b_{n,0}\\
\vdots&\vdots&\vdots&\ddots&\vdots\\
0&0&b_{3,0}&\ldots& b_{n,n-3}\\
0&b_{2,0}&b_{3,1}&\ldots& b_{n,n-2}\\
b_{1,0}&b_{2,1}&b_{3,2}&\ldots& b_{n,n-1}\\
\end{array}\right)$$
with entries from the field  $K$, ~  $b_{j,0}\ne 0$ for $1\leq j\leq
n$. It is well known  that  $\Lc\subset \Omega$.  Well known that
the subset $\Omega$ is invariant with respect to  $\Ad_U$ and each
$\Ad_U$-orbit of  $\Omega $ intersects  $\Lc$  in a unique point.

 Since  $\Omega$ is dense in  $M$, we see that any  $f\in K(M)^U$ is uniquely determined by its restriction on
 $\Lc$. The map
 $\pi$,  that takes  a  rational function  $f\in K(M)^U$ its restriction on $\Lc$,
 is an embedding
 of the field
 $K(M)^U$ into the field $K(\Lc)$.

Construct the matrix
$$ \Sb =
\left(\begin{array}{ccccc}
0&0&0&\ldots&s_{n,0}\\
\vdots&\vdots&\vdots&\ddots&\vdots\\
0&0&s_{3,0}&\ldots& s_{n,n-3}\\
0&s_{2,0}&s_{3,1}&\ldots& s_{n,n-2}\\
s_{1,0}&s_{2,1}&s_{3,2}&\ldots& s_{n,n-1}\\
\end{array}\right),$$
where $s_{k,i}$  is a restriction of the corresponding element of
the matrix  $\Xb$ on $\Lc$ (more precisely,  $x_{a,b}$ with
$a=i-k+n+1$,\, $b=k$).
 It is obvious that  $s_{k,i}(B) = b_{k,i}$.
The system of coordinate functions  $\{s_{k,i}\}$ is algebraically
independent and generate the field  $K(\Lc)$.

 Order    $\{ s_{i,j}\}$ in the following way
$$ s_{1,0}<\ldots <s_{n,0}<s_{2,1}< s_{3,2}<s_{3,1}<
\ldots <s_{n,n-1}<\ldots < s_{n,1}.\eqno (2)$$ {\bf Proposition 2}.
For any pair $(k,i)$, where $1\leq k\leq n,\, 0\leq i\leq k-1$ there
exist the rational functions  $\phi_{k,i}\ne 0$ and $\psi_{k,i}$ of
coordinate functions   $\{s_{a,b}\}$ that are smaller than $s_{k,i}$
in the sense  (2) such that $$\pi(J_{k,i}) = \phi_{k,i} s_{k,i} +
\psi_{k,i}.\eqno (3)$$  {\bf Proof}. We use the induction method on
the coordinate functions  given by  (2). The statement is true for
$i=0$. Assume that the statement is proved for all coordinate
functions  $< s_{k,i}$, where $i\geq 1$; let us prove the statement
for $s_{k,i}$.

Denote by  $S_{k,i}$ the block of order $(k-i)$ of the matrix $\Sb$,
that is an intersection of last  $k-i$ rows and columns
$\{i+1,\ldots, k\}$.  Let $C_{i}$ be the left lower corner block of
order
 $i$ of the adjugate matrix  $\Sb^*$. By direct calculations we obtain
$$\pi(J_{k,i}) =
\left|\begin{array}{cc} *&S_{k,i}\\C_{i}&0\\
\end{array}\right|. \eqno(4)$$
The coordinate function  $ s_{k,i}$  is in the right upper corner of
block $S_{k,i}$. All other coordinate functions of $S_{k,i}$ are
smaller than  $ s_{k,i}$. The minor $|S_{k,i}|$ does not equal to
zero. The minor $|C_i|$ is a monomial of $s_{1,0}\cdots s_{n,0}$.
Calculating (4) we obtain (3). $\Box$
\\
EXAMPLE 1. Case n=2,\, \, $\Sb=\left(\begin{array}{cc}
0&s_{2,0}\\
s_{1,0}&s_{21}\end{array}\right)$,\quad
$\Sb^*=\left(\begin{array}{cc}
*&-s_{2,0}\\
-s_{1,0}& 0\end{array}\right)$, \\
$\pi(J_{1,0})= s_{1,0}$,\, $\pi(J_{2,0})= -s_{1,0}s_{2,0}$,\,
$$\pi(J_{2,1}) = \left|\begin{array}{cc} s_{1,0}&s_{2,1}\\
-s_{1,0}&0 \end{array}\right| = s_{1,0}s_{2,1}.$$
\\
EXAMPLE 2. Case n=3,
 $$\Sb=\left(\begin{array}{ccc}
0&0&s_{3,0}\\
0&s_{2,0}&s_{3,1}\\
s_{1,0}&s_{2,1}&s_{3,2}\\
\end{array}\right),\quad
\Sb^*=\left(\begin{array}{ccc}
*&*&s^*_{3,0}\\
*&s^*_{2,0}& 0\\
s^*_{1,0}&0&0\end{array}\right),$$
\\
$$\pi(J_{1,0})= s_{1,0},\quad  \pi(J_{2,0})= -s_{1,0}s_{2,0},\quad
\pi(J_{3,0})= -s_{1,0}s_{2,0} s_{3,0},$$
\\
$$\pi(J_{2,1}) =  - s^*_{1,0}s_{2,1},\quad \quad \pi(J_{3,2}) = \left|\begin{array}{ccc} s_{1,0}&s_{2,1}&s_{3,2} \\
*&s^*_{2,0}&0\\
s^*_{1,0}&0&0 \end{array}\right| = s^*_{1,0}s^*_{2,0}s_{3,2},$$
\\
$$  \pi(J_{3,1}) = \left|\begin{array}{ccc} 0&s_{2,0}&s_{3,1} \\
s_{1,0}&s_{2,1}&s_{3,2}\\
s^*_{1,0}&0&0 \end{array}\right| = s^*_{1,0}\left|\begin{array}{cc}
s_{2,0}s_{3,1}\\
s_{2,1}s_{3,2}\end{array}\right|.$$ {\bf Theorem}. The field of
$U$-invariants of adjoint representation of the group $\GL$ is the
field of rational functions of $\{J_{k,i}:\, 1\leqslant k\leqslant
n,\,\, 0\leqslant i\leq k-1\}$. \\
 {\bf Proof}. As we say above, the
restriction map $\pi$ is an embedding of the field $K(M)^U$ into the
field $K(\Lc)$. It follows from proposition  2 that the system
$\{\pi(J_{k,i}):\, 1\leqslant k\leqslant n,\,\, 0\leqslant i\leq
k-1\}$ is algebraically independent and generate the field $K(\Lc)$.
This implies the statement of  theorem. $\Box$

The authors are grateful  to  M.Brion, C. De Concini, D.F.Timashev,\\
D.A.Shmelkin, E.B.Vinberg  for useful discussions.

K.A.Vyatkina,\\
Samara State University,\\
vjatkina@gmail.com\\
\\
A.N.Panov,\\
Samara State University, \\
apanov@list.ru\\


\begin{thebibliography}{99}

\bibitem{m}
 K.Miyata, \emph{ Invariants of certain groups. 1.}
Nagoya Math. J., 1971, vol. 41,  69-73.
\bibitem{p}
 L.Pukanszky,  Le\c cons sur les repr\'esentations des
groupes, Paris, Dunod, 1967
\end{thebibliography}
\end{document}